\documentclass[12pt]{amsart}
\usepackage[]{amsmath, amsthm, amsfonts, verbatim, amssymb}
\newcommand\C{{\mathbb C}}

\newcommand\dee{\partial}
\renewcommand\O{\Omega}
\newcommand\Obar{\overline{\Omega}}

\numberwithin{equation}{section}

\begin{document}

\title[Proper holomorphic maps]
{The structure of the semigroup of proper holomorphic mappings of a planar domain to the unit disc}
\author[Bell, Kaleem]
{Steven R. Bell and Faisal Kaleem}

\address[]{Mathematics Department, Purdue University, West Lafayette,
IN  47907}
\email{bell@math.purdue.edu}

\address[]{
Department of Mathematics and Physics,
University of Louisiana at Monroe,
Monroe, LA 71209}
\email{kaleem@ulm.edu}

\thanks{Research supported by NSF grant DMS-0305958}  


\begin{abstract}
Given a bounded $n$-connected domain $\O$ in the plane bounded by $n$~non-intersecting
Jordan curves and given one point $b_j$ on each boundary curve, L.~Bieberbach
proved that there exists a proper holomorphic mapping $f$ of $\O$ onto the unit disc
that is an $n$-to-one branched covering with the properties:  $f$ extends
continuously to the boundary and maps each boundary curve one-to-one onto the unit
circle, and $f$ maps each given point $b_j$ on the boundary to the point $1$ in
the unit circle.  We shall modify a proof by H.~Grunsky of Bieberbach's result
to show that there is a rational function of $2n+2$ complex variables that generates
all of these maps.  In fact, we show that there are two Ahlfors maps $f_1$ and $f_2$
associated to the domain such that any such mapping is given by a fixed linear fractional
transformation mapping the right half plane to the unit disc composed with
$c\,R +  i\,C$, where $R$ is a rational function of the $2n+2$ functions
$f_1(z),f_2(z)$ and $f_1(b_1),f_2(b_1),\dots,f_1(b_n),f_2(b_n)$, and
where $c$ and $C$ are arbitrary real constants subject to the condition $c>0$.
We also show how to generate {\it all\/} the proper holomorphic mappings to the unit
disc via the rational function~$R$.
\end{abstract}

\maketitle

\theoremstyle{plain}

\newtheorem {thm}{Theorem}[section]
\newtheorem {lem}[thm]{Lemma}

\hyphenation{bi-hol-o-mor-phic}
\hyphenation{hol-o-mor-phic}

\section{Introduction}
\label{intro}
The Riemann map associated to a simply connected domain $\O\ne\C$
in the complex plane is the conduit for pulling back the lovely
and explicit formulas on the unit disc back to the domain.  Thus,
the Green's function, Poisson kernel, Szeg\"o kernel, and Bergman
kernel can all be expressed very simply and concretely in terms
of a Riemann map.  The line of research in \cite{B1,B2,B3} has given
the Ahlfors mappings associated to a multiply connected domain
in the plane a similar elevated status.  The classical kernel
functions can all be expressed in terms of {\it two\/} Ahlfors
mappings, albeit not as concretely.

An Ahlfors mapping $f$ of a multiply connected domain is an example
of a proper holomorphic mapping to the unit disc, meaning that,
given a compact subset $K$ of the unit disc, $f^{-1}(K)$ must be
a compact subset of the domain.  When the boundary of $\O$ consists of
finitely many non-intersecting Jordan curves, then $f$ extends
continuously to the boundary, and the properness condition is
equivalent to the condition that $f$ maps the boundary of $\O$ into
the boundary of the unit disc.  In the simply connected case, all
the proper holomorphic mappings to the unit disc can be written down
explicitly.  They are all given as finite Blaschke products composed
with a single Riemann map.  The purpose of this paper is to
give a similar explicit description of all possible proper
holomorphic mappings of a multiply connected domain onto the unit
disc.  We first give a technique for generating all the $n$-to-one
proper holomorphic mappings of an $n$-connected domain onto the
unit disc.  These mappings can be thought of as the analogue of
Riemann mappings in the multiply connected setting, since they
are onto and $m$-to-one for the smallest possible $m$.  Afterwards,
we explain how to generate {\it all\/} the proper holomorphic maps
using the basic $n$-to-one maps.

The set of all proper holomorphic mappings of $\O$ to the unit disc
forms a semi-group.  Elements can be multiplied together to get new
proper maps, but division can only take place if the zero set
(counted with multiplicities) of the denominator is a subset of
the zero set of the numerator.

To motivate what follows, we shall take a moment to describe
the set of biholomorphic mappings of the unit disc onto the
right half plane (RHP) that map a given boundary point $b$ to
the point at infinity.  The map $\frac{1+z}{1-z}$, maps the
unit disc biholomorphically onto the RHP and sends the point
$1$ to $\infty$.  Thus, the map
$$\tau_b(z) := \frac{1+\bar b z}{1-\bar b z}$$
is a biholomorphic map from the unit disc to the RHP sending the
boundary point $b$ to $\infty$.  If $\tau(z)$ is any other such
mapping, note that the quotient $\tau(z)/\tau_b(z)$ has a removable
singularity at $b$.  Approach $b$ along the unit circle in
a counterclockwise direction and use the fact that tangent
vectors and inward pointing normals get rotated by an angle via
a conformal map to see that the image under $\tau_b$ of the
point on the circle moves along the imaginary axis and tends to
$-i\infty$ as the point on the circle approaches $b$.
The same reasoning applies to $\tau(z)$.  This
shows that $\tau(z)/\tau_b(z)$ has a real and positive limit as
$z$ tends to $b$.  Hence, the quotient of the residues at $b$
is real and positive, and consequently, there is a positive
constant $c$ such that $\tau(z)-c\tau_b(z)$
has a removable singularity at $b$.  Now $\tau(z)-c\tau_b(z)$ is a holomorphic
function on a neighborhood of the closed unit disc with no
poles on the closed disc which is pure imaginary valued
on the unit circle.  Such a function must be a pure imaginary
constant.  Therefore, $\tau(z)=c\tau_b(z) + iC$ where $c$ is a positive
constant and $C$ is a real constant.  Notice that the function
$$\tau_b(z) = \frac{1+\bar b z}{1-\bar b z} = \frac{b+ z}{b- z}$$
that generates all such maps has a meromorphic extension to the
disc in the $b$ variable.

Now if $\Omega$ is a simply connected domain bounded by a Jordan curve and
$f$ is a Riemann map from $\Omega$ to the unit disc sending a boundary point
$\beta$ to $b$, let
\begin{equation*}
\Phi_\beta(z) =
\frac{f(\beta)+f(z)}
{f(\beta)-f(z)}.
\end{equation*}
All the biholomorphic maps from $\Omega$ to
the right half plane sending $\beta$ to infinity are given by
$$c\,\Phi_\beta(z)+iC$$
where $c>0$ and $C$ is real.  Note that these mappings extend in $z$ and
$\beta$ as rational functions of $f(z)$ and $f(\beta)$.  Our main goal in
this work will be to prove an analogous result in the multiply connected
case where {\it two\/} Ahlfors mappings shall take the place of a single
Riemann map.

Suppose $\O$ is a bounded domain in the plane bounded by $n$
non-intersecting Jordan curves $\gamma_1$, $\gamma_2$,\dots,$\gamma_n$,
and suppose $f$ is a proper holomorphic mapping of $\O$ onto the unit
disc.  It is known that $f$ must extend continuously to the boundary of $\O$
and that for each boundary curve $\gamma_j$, there is a positive
integer $m_j$ such that the extension maps $\gamma_j$ to the unit
circle as an $m_j$-to-one covering map (see \cite{R}).

Given a point $a$ in $\O$, the Ahlfors mapping associated to $a$ is
the holomorphic function $f_a$ mapping $\O$ into the unit disc such
that $f_a'(a)$ is real and maximal.  Ahlfors proved that the Ahlfors
map is a proper holomorphic mapping onto the unit disc which is an
$n$-to-one branched covering map, which extends continuously to the
boundary, and which maps each boundary curve one-to-one onto the
unit circle.  Ahlfors maps extend meromorphically to the double of
$\O$ by virtue of the fact that $f_a(z)=1/\overline{f_a(z)}$ when
$z\in b\O$.  We shall later use the result proved in \cite{B3} that
there are two points $a_1$ and $a_2$ in $\O$ such that the Ahlfors
maps associated to $a_1$ and $a_2$, when extended meromorphically
to the double of $\O$, generate the field of meromorphic functions
on the double.

Given one point $b_j$ from each boundary curve $\gamma_j$,
Bieberbach \cite{Bieb} proved that there exists a proper
holomorphic mapping $f$ which is an $n$-to-one branched covering
map of $\O$ onto the unit disc, and which maps each boundary curve
one-to-one onto the unit circle in such a way that each point $b_j$
gets mapped to $1$.  We will rework a proof of this result given by
H.~Grunsky in \cite{G1} (see also \cite{G2}) to get more information
about the structure of such maps.  It will be convenient to use the
right half plane as the target domain instead of the unit disc.
Indeed, the two settings are completely equivalent because the two
targets are biholomorphic via a simple linear fractional mapping.
Therefore, we shall concern ourselves with constructing a proper
holomorphic mapping $F$ of $\O$ onto the RHP which maps each point
$b_j$ to the point at infinity.  The mapping $F$ is close to being
unique.  Grunsky demonstrates that $F$ is uniquely determined up
to multiplication by a positive constant and addition of an imaginary 
constant.

Let $b=(b_1,\dots,b_n)$ denote the vector of boundary points
and let $F_b$ denote the mapping to the RHP described above.  We
shall call any such mapping a {\it Grunsky map}, since we shall modify
Grunsky's construction to prove our main results.  Our main
theorem yields a rational function which generates all such maps.

\begin{thm}
\label{thm1}
If $\O$ is an $n$-connected domain in the plane bounded by $n$
non-intersecting Jordan curves, there exist two points $a_1$ and
$a_2$ in $\O$ such that the Ahlfors maps $f_1$ and $f_2$ associated
to $a_1$ and $a_2$, respectively, generate the Grunsky maps in the
sense that there is a rational function of $2n+2$ variables such
that any Grunsky map $F_b$ is given by
$$F_b(z)=cR(f_1(z),f_2(z),f_1(b_1),f_2(b_1),\dots,
f_1(b_n),f_2(b_n))+ iC,$$
where $c>0$ and $C$ is real.
\end{thm}

We emphasize that the rational function in the statement of the
theorem is fixed.  Only the constants $c$ and $C$ vary.  By
composing the Grunsky maps in the theorem with the linear fractional
transformation $(z-1)/(z+1)$, a formula is obtained for all the
proper holomorphic $n$-to-one mappings of the domain $\O$ onto the
unit disc.  Each of the points $b_j$ in the boundary curves get
mapped to~$1$.

The set of all proper holomorphic mappings of the domain $\O$
to the right half plane forms a semi-group under addition.  In
section~\ref{sec5} of this paper, we show how to exploit this
feature, together with the rational function of Theorem~\ref{thm1},
to generate all the proper holomorphic mapping to the RHP.

The main results of this paper grew out of the Purdue PhD thesis
\cite{K} of the second author under the direction of the first author.

In the next section, we will outline Grunsky's construction of
the proper mappings so that we can modify it to prove Theorem~\ref{thm1}
in the section that follows.  In section~\ref{sec5}, we explain
how to generate all the proper holomorphic mappings to the RHP.
In the last section of the paper, we mention some applications of
the main results and some avenues for future research.

\section{Grunsky's results}
\label{sec2}
Let $\Omega$ be a bounded domain bounded by $n$ non-intersecting {\it real
analytic\/} Jordan curves, and let $g(z,w)$
denote the classical Green's function associated to $\O$ (with singularity
$-\ln|z-w|$).  The Poisson
kernel associated to $\O$ is given by
$$p(z,w)=\frac{1}{2\pi}\frac{\partial}{\partial n_w}g(z,w)=
-\frac{i}{\pi}\frac{\partial}{\partial w}g(z,w)T(w)$$
where $T(w)$ is the complex number of unit modulus pointing in the
tangential direction to the boundary at $w$ pointing in the direction
of the standard orientation of the boundary.

For a vector $b=(b_1,\dots,b_n)$ of boundary points, Grunsky \cite{G2}
constructed the map $F_b$ by taking a linear combination of the Poisson
kernels $p(z,b_j)$ in such a way that the resulting harmonic function has
a single valued harmonic conjugate.  Grunsky then showed that the
resulting holomorphic function satisfies all the desired requirements.

We shall always denote the boundary curves of $\Omega$ by $\gamma_j$, $j=1,\dots,n$,
and we follow the convention of letting $\gamma_n$ denote the outermost boundary
curve (the curve bounding the unbounded component of the complement of $\Omega$).
Let $C_j$ denote a smooth simple closed curve that is homotopic to $\gamma_j$ which
is obtained by moving in a short distance along the inward pointing normal vector
to the boundary as $\gamma_j$ is traversed in the standard sense.  If $u$ is
harmonic on $\O$, let ${\mathcal P}_j(u)$ denote the increase of a harmonic
conjugate of $u$ along the cycle $C_j$, i.e., the period of $u$ around $C_j$.
Thus
$${\mathcal P}_j(u)=\int_{C_j}\frac{\partial u}{\partial n}ds,$$
where $n$ is the outward normal and $C_j$ is traversed in the standard sense.

It is an elementary consequence of Green's theorem that,
if $u(z)$ is a harmonic function in $\Omega$, then
$$\sum_{j=1}^{n}{\mathcal P}_j(u)=0.$$
Furthermore, $u$ has a single valued harmonic conjugate if and only if
${\mathcal P}_j(u)=0$ for $j=1,\dots,n-1$.  (See \cite[p.~62]{G2} for a proof.)

Grunsky's construction is based on the following elementary linear algebra result,
which is easily proved by induction.

\begin{lem}
\label{lem2.3}
A system of linear equations
$$\sum_{j=1}^{n}c_{ij}x_j=B_i, \qquad i=1,\dots,n$$
such that $c_{ij}<0$ for $i\neq j$, $\sum_{i=1}^{n}c_{ij}>0$ for each $j$, and
$B_i>0, i=1,\dots,n$, has a unique solution.  Furthermore,
the solution satisfies $x_j>0$ for all $x_j$.
\end{lem}

We are now in a position to describe Grunsky's proof of Bieberbach's theorem.

\begin{thm}
\label{thm2.4}
Choose one point $b_j$ in each boundary curve of a bounded $n$-connected domain $\O$
bounded by $n$ non-intersecting Jordan curves.
There exists an $n$-to-$1$ proper holomorphic mapping $F$ from $\Omega$ to the right half
plane, with the usual counting of multiplicities, such that $F$ extends continuously
up to the boundary and $F(b_j)=\infty$ for each $j$.
The function $F$ is unique up to a positive multiplicative and an
imaginary additive constant.
\end{thm}

\begin{proof}
We may assume that the boundary curves of $\O$ are real analytic since domains of
the kind mentioned in the theorem are biholomorphic to domains with real analytic
boundary via a biholomorphic mapping which extends continuously to the boundary.
Define a positive harmonic function in $\Omega$ via
$$u(z):=\sum_{j=1}^{n}a_j p(z,b_j),$$
where $a_j$ are positive numbers to be determined soon and $p(z,w)$ is the Poisson
kernel.  Grunsky showed that the $a_j$ can
be chosen so that $u$ has a single valued harmonic conjugate.
Set, for short, $$p_j(z)=p(z,b_j).$$
Notice that if $i\ne j$, then the
period $\lambda_{ij} := {\mathcal P}_i(p_j)$ satisfies
$$\lambda_{ij}=\int_{\gamma_i}\frac{\partial p_j}{\partial n}ds,$$
since $p_j(z)$ is smooth up to $\gamma_i$ and so $C_i$ can be deformed to $\gamma_i$
in the definition of ${\mathcal P}_i$.
For $i=j$, we use the fact that the sum of all the periods of a harmonic function
is zero to see that
$$\lambda_{jj}:={\mathcal P}_j(p_j) = -\sum_{i=1,i\neq j}^{n} \lambda_{ij}.$$
Notice that this last identity shows that
$$\sum_{i=1}^{n-1}\lambda_{ij}=-\lambda_{nj}.$$

For $u$ to have a single valued harmonic conjugate, it must happen that
${\mathcal P}_i(u)=0$ for $i=1,....,n-1$ (and consequently, that
${\mathcal P}_n(u)=0$ too).  Set $a_n=1$.   In order to make all the periods of $u$
vanish, the coefficients $a_j$ must satisfy
\begin{equation}
\label{eqn2.1}
\sum_{j=1}^{n-1}\lambda_{ij}a_j=-\lambda_{in},\qquad i=1,\dots,n-1.
\end{equation}

We shall next show that the coefficients of this system satisfy the hypothesis of
Lemma~\ref{lem2.3}, and so the system has a unique solution
$a_1,\dots,a_{n-1}$ with each $a_j$ positive.

Since $p_j(z)>0$ for $z\in \Omega$ , and $p_j(z)=0$ for $z\in \gamma_i$ if  $i\neq j$,
the Hopf Lemma implies that $\frac{\partial p_j}{\partial n}(z)<0$ for $z\in \gamma_i$,
$i\neq j$.  Thus, since $\lambda_{ij}$ is an integral of a strictly negative function,
we obtain that $\lambda_{ij}<0$ for $i\neq j$.  Notice that, consequently, $-\lambda_{in}>0$
for $i=1,\dots,n-1$.  Furthermore, $\sum_{i=1}^{n-1}\lambda_{ij}=-\lambda_{nj}$ is positive
for $j=1,\dots,n-1$.
Thus, using Lemma~\ref{lem2.3} with $n-1$ in place of $n$ and taking $\lambda_{ij}$ as the
coefficients $c_{ij}$, $i,j=1,\dots,n-1$, and taking $-\lambda_{in}$ as $B_i$, $i=1,\dots,n-1$,
we get unique positive numbers $a_1,\dots,a_{n-1}$
which satisfy system~(\ref{eqn2.1}).  The resulting harmonic function $u$ has
a single-valued conjugate $v$ on $\O$.  Now the function
$F=u +iv$ is holomorphic in $\Omega$ and for $w\in b\Omega$, $\text{Re\,}F(z)$
approaches $0$ as $z$ approaches $w$ with the exception of one point $b_j$ on each
boundary component where $F(z)$ approaches $\infty$ as $z$ approaches $b_j$.  It is
now a standard matter to see that $F$ is an $n$-to-one proper holomorphic mapping
to the RHP.  The only
parameters in the construction that we can vary are the real positive constant $a_n$ and
the choice of an imaginary constant in the construction of the harmonic conjugate $v$.
This completes the overview of Grunsky's proof.
\end{proof}

\section{Closer scrutiny of Grunsky's coefficients}
\label{sec3}
In this section, we will start by deriving an important relationship between
the periods $\lambda_{ij}$ of the previous section and the classical functions
$F_j'$ which are defined via
$$F_j'=2\frac{\partial \omega_j}{\partial z},$$
where $\omega_j$ is the harmonic measure function which is harmonic in $\Omega$,
has boundary values one on $\gamma_j$ and zero on the other boundary curves.
For a point $z$ in the boundary of $\Omega$, let $T(z)$ denote the complex number
of unit modulus pointing in the tangential direction given by the standard sense
of the boundary.  We shall show that
$$\lambda_{ij}=-iF_i'(b_j)T(b_j).$$
(Like every complex analyst, we view the symbol $i$ to represent a positive
integer as a subscript and the famous complex number elsewhere.)

It is an elementary fact that if $u$ is a real valued harmonic function
that is smooth up to the boundary and constant on the boundary, then
$$\frac{\dee u}{\dee n}\,ds=-2i\frac{\dee u}{\dee z}\, dz=
2i\frac{\dee u}{\dee\bar z}\, d\bar z \qquad\text{on }b\O,$$
where $ds$ denotes arc length.
If we divide this identity by $ds$, we obtain
$$\frac{\dee u}{\dee n}=-2i\frac{\dee u}{\dee z}\, T(z)=
2i\frac{\dee u}{\dee\bar z}\, \overline{T(z)} \qquad\text{on }b\O.$$
Thus, for example,
$$\frac{\dee \omega_j}{\dee n}\,ds=-iF_j'(z)\, dz=
i\overline{F_j'(z)}\, d\bar z \qquad\text{on }b\O,$$
and it follows that
\begin{equation}
\label{eqn3F}
F_j'(z)\, T(z)= -\overline{F_j'(z)}\,\overline{T(z)} \qquad\text{on }b\O.
\end{equation}

Notice that
$$\omega_i(z)=\int_{w\in\gamma_i}p(z,w)ds=
\frac{1}{2\pi}\int_{w\in\gamma_i}\frac{\partial g}{\partial {n_w}}(z,w)ds,$$
and since $g(z,w)$ is harmonic in $w$ and zero on the boundary for for fixed $z$,
it follows that
\begin{equation*}
\omega_i(z)=\frac{i}{\pi}\int_{w\in\gamma_i}\frac{\partial g}{\partial {\bar w}}(z,w)d\bar w
\end{equation*}
We next differentiate this identity with respect to $z$ and multiply by $2$ to obtain the
well known identity
\begin{equation}
\label{eqn3.A}
F_i'(z)=\frac{2i}{\pi}\int_{w\in\gamma_i}\frac{\partial^2}{\partial z\partial\bar w}g(z,w)d\bar w.
\end{equation}

Since each $p_j$ is real valued harmonic and zero on boundary curves $\gamma_i$ with
$i\ne j$, we have that if $i\neq j$, then
\begin{equation}
\label{eqn3.8}
\lambda_{ij}=\int_{\gamma_i}\frac{\partial p_j}{\partial n}ds
=2i\int_{w\in\gamma_i}\frac{\partial p_j}{\partial \bar w}\,d\bar w.
\end{equation}

Recall that $p(w,\zeta)=\frac{1}{2\pi}\frac{\partial}{\partial n_\zeta}g(w,\zeta)
= - \frac{i}{\pi}\frac{\partial}{\partial \zeta}g(z,\zeta)T(\zeta)$.  Hence
$$p_j(w)=p(w,b_j)= - \frac{i}{\pi}\frac{\partial}{\partial z}g(w,b_j)T(b_j)$$
where $z$ denotes the second variable.  So, for $i\ne j$, (\ref{eqn3.8}) gives
\begin{equation*}
\lambda_{ij}=
\frac{2}{\pi}T(b_j)\int_{w\in\gamma_i}\frac{\partial^2}{\partial\bar w\partial z}g(w,b_j)\ d\bar w.
\end{equation*}
The Green's function is symmetric, i.e., $g(z,w)=g(w,z)$ for $(z,w)\in\Obar\times\Obar$ with $z\ne w$,
and so
\begin{equation*}
\lambda_{ij}=
\frac{2}{\pi}T(b_j)\int_{w\in\gamma_i}\frac{\partial^2}{\partial\bar w\partial z}g(b_j,w)\ d\bar w.
\end{equation*}
We may compare this to equation~(\ref{eqn3.A}) to obtain
$$\lambda_{ij}= -i F_i'(b_j)T(b_j).$$

For the case $i=j$, since $\lambda_{jj} = -\sum_{i=1,i\neq j}^{n} \lambda_{ij},$ we obtain
$$
\lambda_{jj}= iT(b_j)\sum_{i=1,i\neq j}^{n}F_i'(b_j)
$$
But $\sum_{i=1}^n \omega_i\equiv 1$, and so $\sum_{i=1}^n F_i'=0$.  Hence, the last equation
reduces to
$$\lambda_{jj}= -i F_j'(b_j)T(b_j),$$
which agrees with our general formula, and the identity is proved.

We next wish to show that the coefficients $a_j$ in the definition of $F$ have extension properties
in the variables $b_j$.  We shall need to use the following elementary fact.  If $h$ and $H$ are
meromorphic functions on $\O$ which extend meromorphically past the boundary of $\O$ without poles
on the boundary, and if $h(z)=\overline{H(z)}$ for $z$ in the boundary, then $h$ extends to the
double of $\O$ as a meromorphic function.  The Hopf Lemma reveals that the functions $F_j'$ do not
vanish on the boundary (since they are non-zero multiples of the normal derivatives of the harmonic measure
functions $\omega_j$ there).  Furthermore, the identity~(\ref{eqn3F}) allows us to note that
$F_j'(z)/F_1'(z)$ is equal to the conjugate of $F_j'(z)/F_1'(z)$ on the boundary.  Hence, it
follows that $F_j'(z)/F_1'(z)$ extends meromorphically to the double of $\O$.  It was proved
in \cite{B3} that there exist two points in $\O$ such that the meromorphic extensions of the
Ahlfors maps $f_1$ and $f_2$ associated to the two points to the double of $\O$ generate the
field of meromorphic functions on the double.  Thus $F_j'(z)/F_1'(z)$ is a rational combination
of $f_1$ and $f_2$.

Consider the system~(\ref{eqn2.1}) of $n-1$ equations in $n-1$ unknowns.  We have already
noted that this system has a unique solution.  Hence the determinant of the matrix of
coefficients $A=\det [\lambda_{ij}]$ is nonzero.  We shall use Cramer's rule to express
$a_j$ in terms of functions with extension properties.  Let $A_j$ denote the matrix obtained
from $A$ by replacing the $j$-th column with the column vector
$(-\lambda_{1,n},\dots,-\lambda_{n-1,n})$.  Now
\begin{equation}
\label{eqnCramer}
a_j=\frac{det(A_j)}{det(A)}.
\end{equation}
Keep in mind that $\lambda_{ij}=-iF_i'(b_j)T(b_j)$.  Let $m$ be a positive integer with
$1\le m\le n-1$ and $m\ne j$.
We may divide columns $m$ of $A_j$ and $A$ by $-iF_1'(b_m)T(b_m)$ for each such $m$ without changing the
value of $a_j$.  In this way, the factors $T(b_m)$ cancel out from the columns and we are left
with column entries
of the form $F_i'(b_m)/F_1'(b_m)$, which are rational functions of
$f_1(b_m)$ and $f_2(b_m)$.  We must treat the $j$-th columns differently.  Multiply the
determinant $A_j$ by unity in the form $F_1'(b_n)T(b_n)/[F_1'(b_n)T(b_n)]$ and factor the
denominator into the determinant.
In this way, the factors $T(b_n)$ cancel out from the column and we are left with column entries
of the form $iF_i'(b_n)/F_1'(b_n)$, which are rational functions of
$f_1(b_n)$ and $f_2(b_n)$.  Now multiply the
determinant $A$ by unity in the form $F_1'(b_j)T(b_j)/[F_1'(b_j)T(b_j)]$ and factor the
denominator into the determinant.
In this way, the factors $T(b_j)$ cancel out from the column and we are left with column entries
of the form $-iF_i'(b_j)/F_1'(b_j)$, which are rational functions of
$f_1(b_j)$ and $f_2(b_j)$.  We obtain that
\begin{equation}
\label{eqn3aj}
a_j=\frac{F_1'(b_n)T(b_n)}{F_1'(b_j)T(b_j)} Q_j(b_1,\dots,b_n)
\end{equation}
where $Q_j(b_1,\dots,b_n)$ is a rational function of $f_1(b_k)$ and $f_2(b_k)$ for $k=1,2,\dots,n$.
In particular, $Q_j$ extends meromorphically to the double of $\O$ in each $b_k$ separately.

\section{Proof of the main theorems}
\label{sec4}
We continue to assume that $\O$ is a bounded domain in the plane bounded
by $n$ non-intersecting real analytic curves.
We shall need to use a rather long formula proved in \cite{B2} that
relates the Poisson kernel to the Szeg\H o kernel $S(z,w)$ and the Garabedian
kernel $L(z,w)$.  Before we write the formula, we remark that we shall need
to use the following facts about the Szeg\H o and Garabedian kernels on a
domain with real analytic boundary (proofs of which can be found in \cite{B1}).
The kernel $S(z,w)$ extends holomorphically past the boundary in $z$ for each
fixed $w$ in $\O$.  It extends meromorphically past the boundary in $z$ for
each fixed $w$ in $b\O$; in fact, it extends holomorphically past $b\O-\{w\}$
and has only a simple pole at the point $w$.  Furthermore $S(z,w)\ne0$ if
$z\in b\O$ and $w\in\O$.  If $w\in b\O$, then $S(z,w)$ has exactly $n-1$
simple zeroes, one on each boundary curve different from
the one containing the point $w$.  The kernel $L(z,w)$ has a simple pole in
$z$ at the point $w\in\O$.  It extends holomorphically past the boundary in $z$
for each fixed $w$ in $\O$.  It extends meromorphically past the boundary in $z$
for each fixed $w$ in $b\O$; in fact, it extends holomorphically past $b\O-\{w\}$
and has only a simple pole at the point $w$.  Furthermore $L(z,w)\ne0$ if
$z,w\in \O$ with $z\ne w$.  If $w\in b\O$, then $L(z,w)$ has exactly $n-1$
simple zeroes, one on each boundary curve different from
the one containing the point $w$ (and these zeroes agree with those of the
Szeg\H o kernel).  Finally, $S(z,w)$ is in
$C^\infty$ of $\Obar\times\Obar$ minus the boundary diagonal
$\{(z,z): z\in b\O\}$ and $L(z,w)$ is in
$C^\infty$ of $\Obar\times\Obar$ minus the diagonal
$\{(z,z): z\in \Obar\}$.

It is shown in \cite[p. 1358-1362]{B2} that there is a
point $a$ in $\O$ such that $S(z,a)$ has exactly $n-1$ simple zeroes, and
such that the Poisson kernel is given by
\begin{eqnarray}
\label{eqn4.1}
p(z,w) & = & 2\text{Re }\left[\frac{S(z,w)L(w,a)}{L(z,a)} - \sum_{i=1}^{n-1}\sigma_i(z)F_i'(w)T(w)\right] \\
& & + \frac{|S(w,a)|^2}{S(a,a)} + \sum_{i=1}^{n-1}\tau_i(z)F_i'(w)T(w) \nonumber
\end{eqnarray}
where the functions $\sigma_i$ and $\tau_i$ are functions in $C^\infty(\Obar)$.
Recall that the real part of our proper holomorphic $F$ is given by
\begin{equation}
\label{eqn4.2}
\text{Re }F(z) = \sum_{j=1}^n a_j p(z,b_j),
\end{equation}
where $a_n=1$ and the other coefficients $a_j$ are real and positive and satisfy
\begin{equation*}
\sum_{j=1}^{n-1}\lambda_{ij}a_j=-\lambda_{in},\qquad i=1,\dots,n-1,
\end{equation*}
and where, as shown in \S3, $\lambda_{ij}=-iF_i'(b_j)T(b_j)$.  Inserting these values of
the $\lambda_{ij}$ in the last equation yields
\begin{equation}
\label{eqn4.3}
0 = F_i(b_n)T(b_n) + \sum_{j=1}^{n-1}a_j F_i'(b_j)T(b_j),\qquad i=1,\dots,n-1,
\end{equation}
We shall next combine these results to prove the following theorem.

\begin{thm}
\label{thm3.3}
The function $F$ constructed in the proof of Theorem~\ref{thm2.4} is given by
$$F=\sum_{j=1}^{n}\left(2a_j\frac{S(z,b_j)L(b_j,a)}{L(z,a)}+a_j\frac{|S(b_j,a)|^2}{S(a,a)}\right)+iC,$$
where $a$ is a point in $\Omega$ such that the $n-1$ zeros of $S(z,a)$ are distinct and
simple, and C is a real constant.
\end{thm}

\begin{proof}
When we insert formula~(\ref{eqn4.1}) for the Poisson kernel into equation~(\ref{eqn4.2})
and make note of the vanishing of sums of the form~(\ref{eqn4.3}), and use the fact that
the $a_j$ are real, we obtain
$$\text{Re }F(z)=
\sum_{j=1}^{n}
\left(
2\text{Re }\left[
a_j\frac{S(z,b_j)L(b_j,a)}{L(z,a)}
\right]
+a_j\frac{|S(b_j,a)|^2}{S(a,a)}\right).$$
Notice how all the indeterminate functions $\sigma_i$ and $\tau_i$ have disappeared
conveniently!  Notice also that both sides of the last equation are equal to the
real parts of holomorphic functions.  Therefore, the holomorphic functions differ
by an imaginary constant, and the theorem is proved.
\end{proof}

We remark here that a shorter formula for the mapping in Theorem~\ref{thm3.3}
can be obtained by letting the point $a$ approach a point $a_0$ in the boundary.
Indeed, the set of points $a\in\Obar$ where the zeroes of $S(z,a)$ may not be
distinct and simple is a finite subset of $\O$ (see \cite{B2}).  Note that
$L(z,a)$ has a simple pole at $a$ when $a\in\Obar$, and therefore the constant
$C$ is a function of $a$ given by
$$iC(a)=F(a)- \sum_{j=1}^{n} a_j\frac{|S(b_j,a)|^2}{S(a,a)}.$$
Now $S(a,a)$ tends to infinity as $a$ approaches the boundary.  Furthermore,
$L(z,a)$ is non-zero if $z\in\Obar$ and $a\in\O$ with $z\ne a$, and when $a_0$ is
on the boundary, $L(z,a_0)$ has exactly $n-1$ zeroes, one on each boundary curve
different from the curve containing $a_0$ (see \cite{B2}).  In particular
$L(b_j,a_0)$ is non-zero if $a_0$ is on $\gamma_j$ and $a_0\ne b_j$.  As we let
$a$ approach a boundary point $a_0$ in $\gamma_j$ different from $b_j$, we obtain
that the mapping $F$ is given by
$$F=\sum_{j=1}^{n} 2a_j\frac{S(z,b_j)L(b_j,a_0)}{L(z,a_0)} +iC$$
where $C=\lim_{a\to a_0} C(a) =-iF(a_0)$.  An interesting byproduct of the proof
is that the sum maps $a_0$ to zero.  Hence we have proved the following theorem.

\begin{thm}
\label{thm3.3a}
The function given by
$$\sum_{j=1}^{n} 2a_j\frac{S(z,b_j)L(b_j,a_0)}{L(z,a_0)},$$
where the points $b_j$ are boundary points in $\gamma_j$,
$j=1,\dots,n$, and $a_0$ is a point in the boundary different
from any of the $b_j$, is a Grunsky map that maps each $b_j$
to infinity and the point $a_0$ to zero.
\end{thm}

Another formula for the Poisson kernel proved in \cite{B2},
\begin{equation}
\label{eqn4.2a}
p(z,w) = \frac{|S(z,w)|^2}{S(z,z)} + \sum_{i=1}^{n-1}\nu_i(z)F_i'(w)T(w),
\end{equation}
where the functions $\nu_i$ are in $C^\infty(\Obar)$, leads to another interesting
relationship between $F$ and the Szeg\H o kernel.  Indeed, repeating the proof of
Theorem~\ref{thm3.3} with this expression for the Poisson kernel in place of
(\ref{eqn4.1}) yields the identity
$$\text{Re }F(z)=
\sum_{j=1}^{n}
a_j\frac{|S(z,b_j)|^2}{S(z,z)}.$$
It is interesting to note that, although the individual functions on the right
hand side of this identity are not harmonic, the sum is.

We now turn to completing the proof of Theorem~\ref{thm1}.  The main tools in what
remains of the proof are two formulas proved in \cite{B2} for the Szeg\H o and
Garabedian kernels and a fact about the meromorphic extension of certain types of
functions to the double.  To write down the formulas, recall that the point $a$
used in Theorem~\ref{thm3.3} is such that the $n-1$ zeroes of $S(z,a)$ in the
$z$ variable in $\O$ are distinct and simple.  Denote them by $a_1,a_2,\dots,a_{n-1}$,
and let $a_0$ denote $a$.
The Szeg\H o kernel is given by
\begin{equation}
\label{eqnseg}
S(z,w) = \frac{1}{1-f(z)\overline{f(w)}}\left( c_0 S(z,a)\overline{S(w,a)} +
\sum_{i,j=1}^{n-1} c_{ij}S(z,a_i)\overline{S(w,a_j)}\right),
\end{equation}
where $f$ denotes the Ahlfors map associated to $a$, and the Garabedian kernel
$L(z,w)$ is given by
\begin{equation}
\label{eqngar}
L(z,w)=\frac{f(w)}{f(z)-f(w)}\left( c_0 S(z,a)L(w,a) +
\sum_{i,j=1}^{n-1} \bar c_{ij}S(z,a_i)L(w,a_j)\right).
\end{equation}
The extension fact we shall need is that if $G_j(z)$ and $H_j(z)$ are meromorphic
functions on $\O$ which extend meromorphically to a neighborhood of $\Obar$ such
that
\begin{equation}
\label{eqncond}
G_j(z)T(z)=\overline{H_j(z)}\overline{T(z)}\quad\text{for }z\in b\O,\text{ and }j=1,2,
\end{equation}
then $G_1/G_2$ extends meromorphically to the double of $\O$ (because
$G_1/G_2$ is equal to the conjugate of $H_1/H_2$ on the boundary).  Important
functions on $\O$ that satisfy condition~(\ref{eqncond}) include the functions
$F_j'(z)$, $S(z,a_i)S(z,a_j)$, and $S(z,a_i)L(z,a_j)$.  Indeed, equation~(\ref{eqn3F})
shows that $F_j'$ has the property and the well known identity
$$\overline{S(z,w)}=\frac{1}{i}L(z,w)T(z)\quad\text{for }z\in b\O, w\in\O$$
can be used to see that 
$$S(z,a_i)S(z,a_j)T(z)$$
is equal to the conjugate of
$$-L(z,a_i)L(z,a_j)T(z)$$
on the boundary, and
$$S(z,a_i)L(z,a_j)T(z)$$
is equal to the conjugate of
$$-L(z,a_i)S(z,a_j)T(z)$$
on the boundary.

To continue the proof, we may strip the proper map $F$ in Theorem~\ref{thm3.3} of its
imaginary constant, and we may divide $F$ by the positive constant $iF_1'(b_n)T(b_n)$
(which is equal to minus one half the normal derivative of $\omega_1$ at $b_n$).  We next insert
the values for $a_j$ shown in equation~(\ref{eqn3aj}) into the formula for $F$ to obtain
$$F=\frac{1}{i}\sum_{j=1}^{n}\left(
2Q_j\frac{S(z,b_j)L(b_j,a)}{F_1'(b_j)T(b_j)L(z,a)}+
Q_j\frac{|S(b_j,a)|^2}{F_1'(b_j)T(b_j)S(a,a)}\right).$$
Notice that $S(z,b_j)/T(b_j)=\overline{S(b_j,z)T(b_j)}=\frac{1}{i}L(b_j,z)$, and
so
$$F=-\sum_{j=1}^{n}\left(
2Q_j\frac{L(b_j,z)L(b_j,a)}{F_1'(b_j)L(z,a)}+
Q_j\frac{L(b_j,a)S(b_j,a)}{F_1'(b_j)S(a,a)}\right).$$
We next replace the Szeg\H o and Garabedian kernels in this expression by the expressions
given by~(\ref{eqnseg}) and~(\ref{eqngar}).  The functions $Q_j$ extend
to the double in $z$ and each $b_j$, and the Ahlfors map $f$
extends to the double too.  We shall
factor all those terms involving functions which extend to the double in $z$ or any $b_j$
out in front of each term in the large sum and concentrate on what remains, which are
constants times terms of the form
$$\frac{S(b_j,a_k)L(z,a_i)S(b_j,a_m)}{F_1'(b_j)L(z,a)}\quad\text{ and }\quad
\frac{L(b_j,a)S(b_j,a)}{F_1'(b_j)}.$$
The second expression extends meromorphically to the double in the variable $b_j$
because it is a quotient of functions of $b_j$ satisfying condition~(\ref{eqncond}).  We may
simplify the first expression by noting that quotients of the form $L(z,a_i)/L(z,a)$
extend meromorphically to the double in $z$ since they are equal to the conjugate of
$S(z,a_i)/S(z,a)$ on the boundary.  Hence, we may factor that term out into the front
matter involving functions that extend to the double.  Thus, we are left with terms of
the form
$$\frac{S(b_j,a_k)S(b_j,a_m)}{F_1'(b_j)},$$
and they extend meromorphically to the double in the $b_j$ variable because they
are a quotient of functions satisfying condition~(\ref{eqncond}).  We may finally
conclude that $F$ is a rational combination of functions of $z$ which extend
meromorphically to the double and functions of $b_j$ which extend meromorphically
to the double.  Since, as mentioned before, meromorphic functions on the double
are generated by two Ahlfors maps, $f_1$ and $f_2$, the proof of Theorem~\ref{thm1}
is complete in case $\O$ has real analytic boundary curves.  If $\O$ is
bounded by $n$ non-intersecting Jordan curves, it is a standard device to
map $\O$ biholomorphically to a domain with real analytic boundary (see Grunsky
\cite{G2}).  Since the biholomorphic mapping extends continuously to the
boundary, and since Ahlfors maps composed with biholomorphic maps are
themselves Ahlfors maps, all the results readily carry over to the more
general setting.  This completes the proof.

\section{Proper holomorphic mappings of higher mapping degree}
\label{sec5}
The Grunsky maps can be used to build up proper holomorphic mappings
with higher mapping degrees.  Indeed, suppose $\O$ is a bounded domain
bounded by $n$ non-intersecting smooth real analytic curves and fix
points $b_1,\dots,b_{n-1}$ in the boundary curves
$\gamma_1,\dots,\gamma_{n-1}$ respectively.  Now choose {\it two\/}
distinct points $b_{n,1}$ and $b_{n,2}$ in $\gamma_n$.  Let $F_j$
be a Grunsky map that takes $b_1,\dots,b_{n-1},b_{n,j}$ to the point
at infinity, $j=1,2$.  For any two positive constants $c_1$ and $c_2$,
the mapping $F=c_1F_1+c_2F_2$ is a proper holomorphic mapping of $\O$ to
the right half plane that is an $(n+1)$-to-one branched covering.  It
maps $\gamma_n$ two-to-one onto the imaginary axis union the point at
infinity, and each of the other boundary curves one-to-one onto the
imaginary axis union~$\infty$.  Let ${\mathcal B}=F^{-1}(\infty)$.  We
may now add another boundary point that maps to infinity as follows.
Pick a boundary curve $\gamma_k$ and a point $\beta_k$ on it that is
distinct from all the other boundary points in ${\mathcal B}$.  Now pick 
points $\beta_j$ in ${\mathcal B}\cap\gamma_j$ for $1\le j\le n$ with
$j\ne k$ and let $F_3$ denote the Grunsky map associated to the
boundary points $\beta_1,\dots,\beta_n$.  Now $F+c_3F_3$, where $c_3$ is
any positive constant, is a proper holomorphic mapping of $\O$ onto
the RHP that maps the points in ${\mathcal B}\cup\{\beta_k\}$ to
$\infty$.  By adding a boundary point one at a time in this manner,
we may build up proper holomorphic maps to the RHP with arbitrarily
high mapping degree which are $K_j$-to-one maps of $\gamma_j$ onto
the imaginary axis union $\infty$ where $K_j$ are arbitrary positive
integers.

We shall show in a moment that we might not generate all the proper
holomorphic mappings in this manner.  We might need to allow some of
the coefficients to be negative.

Before we continue, we will take a closer look at the process of
adding a single extra boundary point to the list attached to a
Grunsky map.  As above, fix points $b_1,\dots,b_{n-1}$ in the boundary
curves $\gamma_1,\dots,\gamma_{n-1}$, respectively, and choose two points
$b_{n,1}$ and $b_{n,2}$ in $\gamma_n$.  Let $F_j$ be a Grunsky map that
takes $b_1,\dots,b_{n-1},b_{n,j}$ to the point at infinity, $j=1,2$.
As mentioned above, for any two positive constants $c_1$ and $c_2$,
the mapping $c_1F_1+c_2F_2$ is a proper holomorphic mapping of $\O$ to
the right half plane that is an $(n+1)$-to-one branched covering that
maps the list of points to $\infty$.  We shall now show that all such
maps are given by $c_1F_1+c_2F_2 +iC$ where $c_1$ and $c_2$ are positive
constants and $iC$ is an imaginary constant.  Let ${\mathcal B}$ denote
the set $\{b_j: j=1,\dots,n-1\}\cup\{b_{n,k}: k=1,2\}$.

We now construct another proper
holomorphic map to the RHP that maps the $n+1$ boundary points in
${\mathcal B}$ to the point at infinity using Grunsky's technique.  Indeed, let
$$u(z):=\sum_{j=1}^{n-1}a_j p(z,b_j) +\sum_{k=1}^{2}A_k p(z,b_{n,k}),$$
where we consider the $A_k$ to be arbitrary positive
constants and we determine $a_1,\dots,a_{n-1}$ in order to make the
periods of $u$ vanish.  The required conditions on the coefficients
are
$$\sum_{j=1}^{n-1}a_j F_i'(b_j)T(b_j)= \sum_{k=1}^2-A_k F_i'(b_{n,k})T(b_{n,k}),$$
for $i=1,\dots,n-1$.  This system has a unique solution with each $a_j>0$.
When we use Cramer's Rule as in section~\ref{sec3} and the linearity of
the determinate in the $j$-th column to express the
$a_j$ in terms of functions that extend to the double, we obtain
$$
a_j=
\sum_{k=1}^2 A_k\frac{F_1'(b_{n,k})T(b_{n,k})}{F_1'(b_j)T(b_j)} Q_{j,k}(b_1,\dots,b_{n-1},b_{n,k}),
$$
where the functions $Q_{k,j}$ are rational functions of
$f_1(b_k)$ and $f_2(b_k)$ for $k=1,2,\dots,n$.
We may divide each $A_k$ by the positive number $iF_1'(b_{n,k})T(b_{n,k})$ and maintain
their arbitrary positive nature (since the product is minus the normal derivative of
$\omega_1$ at a boundary point where $\omega_1$ takes its absolute minimum).
The coefficients obtained in this way are identical to the ones obtained in section~\ref{sec4}
for individual Grunsky maps.
Indeed, the proper holomorphic mapping $F$ we obtain in this way is equal, up to the
addition of an imaginary constant, to $A_1$ times the Grunsky map for
$b_1,\dots,b_{n-1},b_{n,1}$ constructed in section~\ref{sec4} plus $A_2$ times the Grunsky map
for $b_1,\dots,b_{n-1},b_{n,2}$.  It is clear that $N$ boundary points $b_{n,k}$,
$k=1,\dots,N$ can be added on $\gamma_n$ under the same procedure to obtain a proper
holomorphic map that maps $\gamma_n$ to the boundary of the RHP in an $N$-to-one manner
and all the other boundary curves in a one-to-one manner.

There was nothing special about choosing the boundary curve $\gamma_n$ above.
Extra points can be added to any boundary curve in the same way.

We shall now explain how to generate {\it all\/} the proper holomorphic
mappings to the RHP.

\begin{thm}
All the proper holomorphic mappings of $\O$ to the RHP are given by
linear combinations of Grunsky maps $F=\sum_{k=1}^N c_k F_k$ subject to the
conditions that all the $c_k$ are real and that when the real part of the
map is decomposed into a sum of the form
$$\text{Re }F(z)=\sum_{j=1}^M a_j p(z,b_j),$$
where the $b_j$ are distinct boundary points, there must be at least one
$b_j$ on each boundary curve, and each coefficient $a_j$ must be strictly
positive.
\end{thm}

We remark that the positivity condition in the theorem is just a
sequence of linear inequalities on the coefficients $c_k$.  Taking each
$c_k$ to be positive always leads to a proper map.  We shall see during
the proof, however, that it is possible to construct proper maps where
some of the $c_k$ might be negative.

To prove the theorem, we use induction on the number of boundary points
in ${\mathcal B}=F^{-1}(\infty)$ where $F$ is a proper holomorphic mapping
of $\O$ onto the RHP.  We handled the case of $N=n$ and $n+1$ boundary
points above (where $n$ is the connectivity of $\O$).  Assume the theorem
is true in the case of $N$ boundary points and suppose $F$ is a proper
holomorphic map to the RHP of mapping degree $N+1$.
We know that ${\mathcal B}$ must contain at least one boundary point from
each boundary curve.  Choose points $b_1,\dots,b_n$ from ${\mathcal B}$ in
$\gamma_1.\dots,\gamma_n$, respectively.  Let $f$ denote the Grunsky map
associated to this sequence of boundary points.
Since $N>n$, there is at least one boundary curve, say $\gamma_k$, that
contains more than one point.  Let
${\mathcal B}_0={\mathcal B}-\{b_k\}$ and let $F_0$ be a proper holomorphic
mapping (as constructed above) such that $F_0^{-1}(\infty)={\mathcal B}_0$.
We now claim that it is possible to choose positive constants $c$ and $c_0$
so that $F+c_0F_0-cf$ is a proper holomorphic mapping to the RHP of mapping
degree~$N$.  We may decompose the real part of $F+c_0F_0$ as
\begin{eqnarray*}
\text{Re }(F+c_0F_0) & = & a_kp(z,b_k)+
\sum_{j=1,j\ne k}^n (a_j +c_0\alpha_j)p(z,b_j)\\
& & +\sum_{m=1}^M (A_m +c_0B_m)p(z,\beta_m),
\end{eqnarray*}
where all the coefficients $a_j$, $\alpha_j$, $A_m$, and $B_m$ are positive.
By choosing $c_0$ sufficiently large (and positive), we may choose a $c>0$
so that the similar decomposition for $F+c_0F_0-cf$ has a zero coefficient
in front of the term $p(z,b_k)$ and positive coefficients in front of the
other $N$ Poisson kernel terms.  Since the pole at $b_k$ is removed, the
resulting map is a proper holomorphic map with mapping degree $N$.  Our
induction hypothesis yields that $F+c_0F_0-cf$ is a linear combination of
Grunsky maps.  Our construction of $F_0$ is also a linear combination of
Grunsky maps.  It follows that $F$ is a linear combination of Grunsky maps.
The positivity condition in the theorem is a necessary feature of any
proper holomorphic map to the RHP.  This completes the proof.

\section{Applications and remarks}
\label{sec6}
It is particularly easy to find ``primitive pairs'' among the Grunsky maps.
Indeed, the Grunsky maps extend meromorphically to the double by simple
reflection.  They extend to be $n$-to-one mappings of the double to the
extended complex plane.  Given a Grunsky map $F_1$, choose points $b_j$, one
in each boundary curve, so that $\{F_1(b_j): j=1,\dots,n\}$ consists of
$n$ distinct points in the finite complex plane.  (This is easy to do since
$F_1$ maps each boundary curve one-to-one onto the imaginary axis union
the point at infinity.)  Now let $F_2$ be the Grunsky map that maps each
$b_j$ to the point at infinity.  Since $F_1$ separates the points of
$F_2^{-1}(\infty)$, it follows that the extensions of $F_1$ and $F_2$ to
the double generate all the meromorphic functions on the double, i.e., they
form a primitive pair (see Farkas and Kra \cite{FK}).

The process of constructing Grunsky maps can be carried out on a finite
Riemann surface and the argument above can be used to construct 
pairs $F_1$ and $F_2$ of Grunsky maps that extend to the double to form
primitive pairs for the double.  When this line of reasoning is combined
with results in \cite{B3}, it follows that the $(1,1)$-form that is the
Bergman kernel on a finite Riemann surface can be expressed via
$$K(z,w)=dF_1(z)R(F_1(z),F_2(z),\overline{F_1(w)},\overline{F_2(w)})d\overline{F_1(w)},$$
where $R$ is a rational function of four complex variables.

Finally, we remark that although we expressed the Grunsky maps in
terms of two Ahlfors maps, we could just as easily have expressed
them in terms of rational combinations of {\it any\/} two meromorphic
functions on $\O$ that extend to the double to form a primitive
pair for the double.  When $\O$ happens to be a quadrature domain,
the function $z$ and the Schwarz function $S(z)$ (which is an algebraic
function that is meromorphic on $\O$) are a particularly appealing
choice (see Gustafsson \cite{G}).  Since any smooth finitely
connected domain is biholomorphic to a quadrature domain that is
$C^\infty$ close by (see \cite{B4,B5}), we conclude that it is possible
to make subtle holomorphic changes of variables so that the Grunsky
maps can be expressed in terms of rational functions of $z$ and
the Schwarz function.  This is very analogous to the prime example
of a quadrature domain, the unit disc, where the Schwarz function
is $S(z)=1/z$ and the Grunsky maps are simple linear fractional
transformations in $z$ and $b$ as shown in the introduction.

\end{document}